\documentclass[12pt,a4paper]{amsart}
\usepackage{mathrsfs}

\newtheorem{theo+}              {Theorem}           [section]
\newtheorem{prop+}  [theo+]     {Proposition}
\newtheorem{coro+}  [theo+]     {Corollary}
\newtheorem{lemm+}  [theo+]     {Lemma}
\newtheorem{exam+}  [theo+]     {Example}
\newtheorem{rema+}  [theo+]     {Remark}
\newtheorem{defi+}  [theo+]     {Definition}

\newenvironment{theorem}{\begin{theo+}}{\end{theo+}}

\newenvironment{corollary}{\begin{coro+}}{\end{coro+}}
\newenvironment{lemma}{\begin{lemm+}}{\end{lemm+}}

\usepackage{amsthm}
\theoremstyle{plain} \theoremstyle{remark}

\newtheorem{example}{Example}

\def\E{/\kern-1.0em \equiv }

\author{Ze-Ping Wang$^{*}$ }
\author{Ye-Lin Ou$^{**}$ }
\address{Department of Mathematics $\&$ Physics,\newline\indent Yunnan
Wenshan University,\newline\indent No. 2 Xuefu Road Wenshan County
Wenshan, Yunnan 653000,\newline\indent People's Republic of China
\newline\indent E-mail:zeping.wang@gmail.com \;(Wang)\\\newline\indent  Department of
Mathematics,\newline\indent Texas A $\&$ M University-Commerce,
\newline\indent Commerce TX 75429,\newline\indent USA.\newline\indent
E-mail:yelin$\_$ou@tamu-commerce.edu \;(Ou)}
\thanks{*Supported by Yunnan
Wenshan University Research Project 09WSY03. The author is also
grateful to the Department of Mathematics, Texas A $\&$ M
University-Commerce for the hospitality he received during a visit in Fall 2009 during which this work was done. \\
\indent** Supported by Texas A $\&$ M University-Commerce ``Faculty
Research Enhancement Project" (2009-10) .}
\date{02/23/2010}
\begin{document}
\title[Biharmonic Riemannian submersions]{Biharmonic Riemannian submersions from $3$-manifolds}

\subjclass{58E20, 53C43} \keywords{Biharmonic maps, Riemannian
submersions, harmonic morphisms, $3$-manifolds.}

\maketitle

\section*{Abstract}
\begin{quote}
{\footnotesize An important theorem about biharmonic submanifolds
proved independently by Chen-Ishikawa \cite{CI} and Jiang \cite{Ji2}
states that an isometric immersion of a surface into $3$-dimensional
Euclidean space is biharmonic if and only if it is harmonic (i.e,
minimal). In a later paper \cite{CMO2}, Cadeo-Monttaldo-Oniciuc
shown that the theorem remains true if the target Euclidean space is
replaced by a $3$-dimensional hyperbolic space form. In this paper,
we prove the dual results for Riemannian submersions, i.e., a
Riemannian submersion from a $3$-dimensional space form of
non-positive curvature into a surface is biharmonic if and only if
it is harmonic.}

\end{quote}

\maketitle

\section{Introduction and the main results}
All manifolds, maps, tensor fields studied in this paper are assumed
to be smooth unless there is an otherwise statement.\\

A {\em biharmonic map} is a map $\varphi:(M, g)\longrightarrow (N,
h)$ between Riemannian manifolds that is  a critical point of the
bienergy
\begin{equation}\nonumber
E^{2}\left(\varphi,\Omega \right)= \frac{1}{2} {\int}_{\Omega}
\left|\tau(\varphi) \right|^{2}{\rm d}x
\end{equation}
for every compact subset $\Omega$ of $M$, where $\tau(\varphi)={\rm
Trace}_{g}\nabla {\rm d} \varphi$ is the tension field of $\varphi$
vanishing of which means the map is harmonic. By computing the first
variation of the functional (see \cite{Ji}) one finds that $\varphi$
is  biharmonic  if and only if its bitension field vanishes
identically, i.e.,
\begin{equation}\label{BT1}
\tau^{2}(\varphi):={\rm
Trace}_{g}(\nabla^{\varphi}\nabla^{\varphi}-\nabla^{\varphi}_{\nabla^{M}})\tau(\varphi)
- {\rm Trace}_{g} R^{N}({\rm d}\varphi, \tau(\varphi)){\rm d}\varphi
=0,
\end{equation}
where $R^{N}$ is the curvature operator of $(N, h)$ defined by
$$R^{N}(X,Y)Z=
[\nabla^{N}_{X},\nabla^{N}_{Y}]Z-\nabla^{N}_{[X,Y]}Z.$$

A submanifold is called a {\bf biharmonic submanifold} if the
isometric immersion that defines the  submanifold is a biharmonic
map. As biharmonic maps include harmonic maps as special cases
biharmonic submanifols generalize the notion of minimal submanifolds
(i.e., minimal isometric immersions). We use {\bf Proper biharmonic
maps (respectively, submanifolds)} to name those biharmonic maps
(respectively, submanifolds) which are not harmonic.\\

A fundamental problem in the study of biharmonic maps is to classify
all proper biharmonic maps between certain model spaces. An example
of this is the following challenging conjecture which is still
open.\\

{\bf Chen's Conjecture \cite{Ch}:} Any biharmonic isometric
immersion $(M^m,
g) \hookrightarrow \mathbb{R}^n$ into Euclidean space is harmonic.\\

Among several cases (see e.g., \cite{Di}, \cite{HV}) that support
the conjecture is the following theorem proved independently by
Chen-Ishikawa \cite{CI} and Jiang \cite{Ji2}.\\

{\bf Theorem.} An isometric immersion $(M^2, g) \hookrightarrow
\mathbb{R}^3$ into Euclidean space is biharmonic if and only if it is harmonic.\\

In a later paper \cite{CMO2}, Cadeo-Monttaldo-Oniciuc shown that the
theorem remains true if the target Euclidean space is replaced by a
$3$-dimensional hyperbolic space form.\\

In this paper, we prove the dual results for Riemannian submersions
and give a complete classification of biharmonic Riemannian
submersions from $3$-dimensional space form. Our main results can be stated as follows.\\

 {\bf Theorem 3.3.} {\em Let
$\pi:(M^3(c),g)\longrightarrow (N^2,h)$ be Riemannian submersion
from a space form of constant sectional curvature $c$. Then,
$\pi$ is biharmonic if and only if it is harmonic}.\\

{\bf Corollary 3.4.} {\em $(1)$ If $\pi:\mathbb{R}^3\longrightarrow
(N^2,h)$ is a biharmonic horizontally homothetic submersion from
Euclidean space, then $(N^2,h)$ is flat and $\pi$ is a composition
of an orthogonal projection $\mathbb{R}^3\longrightarrow
\mathbb{R}^2$
followed by a covering map $\mathbb{R}^2\longrightarrow (N^2,h)$;\\
$(2)$ There exists no biharmonic Riemannian submersion
$\pi:\mathbb{H}^3\longrightarrow (N^2,h)$ no matter what $(N^2,h)$
is}.\\

Applying our results one can easily check the following\\

{\bf Example.} {\em The Riemannian submersion
$\pi:\mathbb{R}^3\longrightarrow (N^2=\mathbb{R}^3/\mathbb{R},h)$
from $\mathbb{R}^3$ onto the orbit space of  a free $1$-parameter
isometric group action of $\mathbb{R}$ on $\mathbb{R}^3$ is not
biharmonic, where the isometric group action is described by
$(s,(z,t))\longrightarrow (e^{is}z,t+s),\;s\in \mathbb{R}, (z,t) \in
\mathbb{C}\times \mathbb{R}\equiv \mathbb{C}\times \mathbb{R}$}. In
fact, if it were, then by our classification results, the Riemannian
submersion would be harmonic and hence a harmonic morphism, and then
a well-known theorem about harmonic morphisms (see, e.g.,
\cite{BW1}) would imply that all fibers of $\pi$ would be geodesics,
which is not the case.

\section{Biharmonic Riemannian submersions from $3$-manifolds}

In this section, we will describe biharmonicity of a Riemannian
submersion from a generic $3$-manifold by using the integrability
data of a special orthonormal frame adapted to a Riemannian
submersion. This is the main tool we use to prove our main theorem.
We also construct a family of proper biharmonic Riemannian
submersions
from $\mathbb{R}^3$ provided with a warped product metric.\\

Let $\pi:( M^3 , g)\longrightarrow (N^2,h)$ be a Riemannian
submersion. A local orthonormal frame is said to be {\bf adapted to
the Riemannian submersion} $\pi$ if the vector fields in the frame
that are tangent to the horizontal distribution are basic (i.e.,
they are $\pi$-related to a local orthonormal frame in the base
space). Such a frame always exists (cf. e.g., \cite{BW1}). Let
$\{e_1,\; e_2, \;e_3\}$ be an orthonormal frame adapted to $\pi$
with $e_3$ being vertical. Then, it is well known (see \cite{On})
that $[e_1,e_3]$ and $ [e_2,e_3]$ are vertical and $[e_1,e_2]$ is
$\pi$-related to $[\varepsilon_1, \varepsilon_2]$, where
$\{\varepsilon_1, \varepsilon_2\}$ is an orthonormal frame in the
base manifold. If we assume that
\begin{equation}\label{RS22}
[\varepsilon_1,\varepsilon_2]=F_1\varepsilon_1+F_2\varepsilon_2,
\end{equation}
for $F_1, F_2\in C^{\infty}(N)$ and use the notations
$f_i=F_i\circ \pi, i=1, 2$. Then, we have
\begin{equation}\label{R1}
\begin{cases}
[e_1,e_3]=\kappa_1e_3,\\
[e_2,e_3]=\kappa_2e_3,\\
[e_1,e_2]=f_1 e_1+f_2e_2-2\sigma e_3.
\end{cases}
\end{equation}
where $\kappa_1,\;\kappa_2\;{\rm and}\; \sigma \in C^{\infty}(M)$.
We will call $ f_1, f_2, \kappa_1,\;\kappa_2\;{\rm and}\; \sigma$
{\bf the integrability data} of the adapted frame of the Riemannian
submersion $\pi$.

\begin{theorem}
Let $\pi:(M^3,g)\longrightarrow (N^2,h)$ be a Riemannian submersion
with the adapted frame $\{e_1,\; e_2, \;e_3\}$ and the integrability
data $ f_1, f_2, \kappa_1,\;\kappa_2\;{\rm and}\; \sigma$. Then, the
Riemannian submersion $\pi$ is biharmonic if and only if
\begin{equation}\label{RSB0}
\begin{cases}
-\Delta^{M}\kappa_1-f_1 e_1(\kappa_2)-e_1(\kappa_2 f_1)-f_2
e_2(\kappa_2)-e_2(\kappa_2 f_2)
\\+\kappa_1\kappa_2 f_1 +\kappa_2^2 f_2
+\kappa_1\{-K^{N}+f_{1}^{2}+f_{2}^{2}\}
=0,\\
-\Delta^{M}\kappa_2+f_1 e_1(\kappa_1)+e_1(\kappa_1 f_1)+f_2
e_2(\kappa_1)+e_2(\kappa_1 f_2)\\
-\kappa_1\kappa_2 f_2-\kappa_1^2f_1
+\kappa_2\{-K^{N}+f_{1}^{2}+f_{2}^{2}\} =0,
\end{cases}
\end{equation}
where
$K^{N}=R^{N}_{1212}\circ\pi=-[e_2(f_1)-e_1(f_2)+f_{1}^{2}+f_{2}^{2}]
$ is the Gauss curvature of Riemannian manifold $(N^2,h)$.
\end{theorem}
\begin{proof}
Let $\nabla$  denote the Levi-Civita connection of the Riemannian
manifold $(M^3,g)$. A straightforward computation using (\ref{R1})
and Koszul formula gives
\begin{eqnarray}\notag
&&\nabla_{e_{1}} e_{1}=-f_1e_2,\;\;\nabla_{e_{1}} e_{2}=f_1
e_1-\sigma e_{3},\;\;\nabla_{e_{1}} e_{3}=\sigma
e_{2},\\\label{Conn} &&\nabla_{e_{2}} e_{1}=-f_2 e_{2}+\sigma
e_3,\;\;\nabla_{e_{2}} e_{2}=f_2 e_{1}, \;\;\nabla_{e_{2}}
e_{3}=-\sigma e_{1},\\\notag && \nabla_{e_{3}}
e_{1}=-\kappa_1e_{3}+\sigma e_{2}, \nabla_{e_{3}} e_{2}= -\sigma
e_{1}-\kappa_2 e_3, \nabla_{e_{3}} e_{3}=\kappa_1 e_{1}+\kappa_2
e_2.
\end{eqnarray}
The tension of the Riemannian  submersion $\pi$ is given by
\begin{equation}\label{RSB6}
\tau(\pi)=\nabla^{\pi}_{e_i}d\pi(e_i)-d\pi(\nabla^{M}_{e_i}e_i)=-d\pi(\nabla^{M}_{e_3}e_3)
=-\kappa_1\varepsilon_1-\kappa_2
\varepsilon_2.
\end{equation}
A straightforward computation using (\ref{Conn}) yields
\begin{equation}\label{RSB7}
\begin{array}{lll}
\sum\limits_{i=1}^2\nabla^{\pi}_{e_i}\nabla^{\pi}_{e_i}\tau(\pi)
=\nabla^{\pi}_{e_i}\nabla^{\pi}_{e_i}(-\kappa_1\varepsilon_1-\kappa_2
\varepsilon_2)\\
=[-e_1e_1(\kappa_1)+\kappa_1f_{1}^2-f_1e_1(\kappa_2)-e_1(\kappa_2
f_1)-e_2e_2(\kappa_1)\\+\kappa_1f_{2}^2-f_2e_2(\kappa_2)-e_2(\kappa_2
f_2)]\varepsilon_1\\
+[-e_1e_1(\kappa_2)+\kappa_2 f_{1}^2+f_1e_1(\kappa_1)+e_1(\kappa_1
f_1)-e_2e_2(\kappa_2)\\+\kappa_2
f_{2}^2+f_2e_2(\kappa_1)+e_2(\kappa_1 f_2)]\varepsilon_2,
\end{array}
\end{equation}

\begin{equation}\label{RSB8}
\begin{array}{lll}
\nabla^{\pi}_{e_3}\nabla^{\pi}_{e_3}\tau(\pi)
=-e_3e_3(\kappa_1)\varepsilon_1 -e_3e_3(\kappa_2)\varepsilon_2,
\end{array}
\end{equation}

\begin{equation}\label{RSB9}
\begin{array}{lll}
\sum\limits_{i=1}^{3}\nabla^{\pi}_{\nabla^{M}_{e_i}e_i}\tau(\pi)\\
=\nabla^{\pi}_{\nabla^{M}_{e_1}e_1}(-\kappa_1\varepsilon_1-\kappa_2
\varepsilon_2)+\nabla^{\pi}_{\nabla^{M}_{e_2}e_2}(-\kappa_1\varepsilon_1-\kappa_2
\varepsilon_2)+\nabla^{\pi}_{\nabla^{M}_{e_3}e_3}(-\kappa_1\varepsilon_1-\kappa_2
\varepsilon_2)\\
=[f_1e_2(\kappa_1)+\kappa_2 f_1f_2 -f_2e_1(\kappa_1)-\kappa_2 f_1f_2
-\kappa_1e_1(\kappa_1)\\-\kappa_1\kappa_2 f_1 -\kappa_2
e_2(\kappa_1)-\kappa_2^2f_2]
\varepsilon_1\\
+[f_1e_2(\kappa_2)-\kappa_1f_1f_2 -f_2e_1(\kappa_2)+\kappa_1f_1f_2
+\kappa_1^2f_1-\kappa_1e_1(\kappa_2) \\+\kappa_1\kappa_2
f_2-\kappa_2 e_2(\kappa_2)] \varepsilon_2,
\end{array}
\end{equation}
and
\begin{equation}\label{RSB12}
\begin{array}{lll}
\sum\limits_{i=1}^3
R^{N}(d\pi(e_i),\tau(\pi))d\pi(e_i)\\
=-\kappa_1\{e_2(f_1)-e_1(f_2)+f_{1}^{2}+f_{2}^{2}\}\varepsilon_1\\
-\kappa_2\{e_2(f_1)-e_1(f_2)+f_{1}^{2}+f_{2}^{2}\}\varepsilon_2.
\end{array}
\end{equation}

Substituting equations (\ref{RSB7})-(\ref{RSB12}) into the bitension
field formula (\ref{BT1}) we obtain
\begin{equation}\label{RSB17}\notag
\begin{array}{lll}
\tau^2(\pi)=\sum\limits_{i=1}^3\{\nabla^{\pi}_{e_i}\nabla^{\pi}_{e_i}\tau(\pi)
-\nabla^{\pi}_{\nabla^{M}_{e_i}e_i}\tau(\pi)
-R^{N}(d\pi(e_i),\tau(\pi))d\pi(e_i)\}\\
=[-\Delta^{M}\kappa_1-f_1e_1(\kappa_2)-e_1(\kappa_2
f_1)-f_2e_2(\kappa_2)-e_2(\kappa_2 f_2)\\+\kappa_1\kappa_2 f_1
+\kappa_2^2f_2 +\kappa_1\{-K^{N}+f_{1}^{2}+f_{2}^{2}\}
]\varepsilon_1\\
+[-\Delta^{M}\kappa_2+f_1e_1(\kappa_1)+e_1(\kappa_1
f_1)+f_2e_2(\kappa_1)+e_2(\kappa_1
f_2)\\
-\kappa_1\kappa_2 f_2-\kappa_1^2f_1
+\kappa_2\{-K^{N}+f_{1}^{2}+f_{2}^{2}\} ]\varepsilon_2,
\end{array}
\end{equation}
from which the theorem follows.
\end{proof}
When the integrability data $\kappa_2=0$ we have the following
corollary which will be used later in the paper.
\begin{corollary}\label{Cor1}
Let $\pi:(M^3,g)\longrightarrow (N^2,h)$ be a Riemannian submersion
with an adapted frame $\{e_1,\; e_2, \;e_3\}$ and the integrability
data $\{ f_1, f_2, \kappa_1,\;\kappa_2,\; \sigma\}$ with
$\kappa_2=0$. Then, the Riemannian submersion $\pi$ is biharmonic if
and only if
\begin{equation}\label{Simp}
\begin{cases}
-\Delta^{M}\kappa_1 +\kappa_1\{-K^{N}+f_{1}^{2}+f_{2}^{2}\}
=0,\\
f_1 e_1(\kappa_1)+e_1(\kappa_1 f_1)+f_2 e_2(\kappa_1)+e_2(\kappa_1
f_2) -\kappa_1^2f_1 =0,
\end{cases}
\end{equation}
\end{corollary}

\begin{example}
For $\varphi(x)=\frac{c_1(1+e^{c_1 x})}{1-e^{c_1 x}}$, and
$\phi(y)=\frac{b_1(1+e^{b_1 y})}{1-e^{b_1 y}}$, and
$\beta(x,y)=ce^{\int\varphi(x)dx+\int\phi(y)dy}$, the Riemannian
submersion
\begin{align}\notag
\pi  : ( \mathbb{R}^2 \times \mathbb{R} , dx^2 +
dy^2+\beta^{-2}(x,y) dz^2) &\to (\mathbb{R}^2 ,dx^2 + dy^2) \\\notag
\phi(x,y,z) =(x,y)
\end{align}
is a proper biharmonic map. In particular, when $\phi(y)=0$, the
example recovers the family of proper biharmonic Riemannian
submersion found in \cite{LO}.
\end{example}

In fact, it is not difficult to check that the orthonormal frame
$\{e_1=\frac{\partial}{\partial x},\;e_2= \frac{\partial}{\partial
y}, \;e_3=\beta \frac{\partial}{\partial z}\}$  on  $( \mathbb{R}^2
\times \mathbb{R} , dx^2 + dy^2+\beta^{-2}(x,y) dz^2)$ is adapted to
the Riemannian submersion $\pi$ with $ d\pi(e_i)=\varepsilon_i, i=1,
2$ and $e_3$ being vertical, where
$\varepsilon_1=\frac{\partial}{\partial x},\;\varepsilon_2=
\frac{\partial}{\partial y},\;$ form an orthonormal frame on the
base space $(\mathbb{R}^2 , dx^2 + dy^2)$. A straightforward
computation gives the Lie brackets
\begin{align}\notag
[e_1,e_3]=f e_3,\;[e_2,e_3]=g e_3 ,\;\; [e_1,e_2]=0,
\end{align}
where $f=(\ln \beta)_{x}, g=(\ln \beta)_{y}\;$.\\
It follows that the integrability data of the Riemannian submersion
$\pi$ are give by
\begin{equation}\notag
\begin{array}{lll}
f_1=f_2=\sigma=0,\;\kappa_1=f,\;\kappa_2=g.
\end{array}
\end{equation}
Substituting these and the curvature $K^{\mathbb{R}^2}=0$ into
Equation (\ref{RSB0}) we conclude that the Riemannian submersion
$\pi$ is biharmonic if and only if
\begin{equation}\label{Ex44}
\begin{array}{lll}
\Delta^{M} f=0,\\
\Delta^{M} g=0.
\end{array}
\end{equation}
Looking for the special solutions of the form
$\ln\beta=\int\varphi(x)dx+\int\phi(y)dy$ we have $f=\varphi(x)$ and
$g=\phi(y)$. Substituting these into  system (\ref{Ex44}) we obtain
a system of ODE:
\begin{equation}\notag
\begin{array}{lll}
\varphi(x)\varphi'(x)- \varphi''(x)=0,\\
\phi(y)\phi'(y)- \phi''(y)=0.
\end{array}
\end{equation}
which has solutions $\varphi(x)=\frac{c_1(1+e^{c_1 x})}{1-e^{c_1
x}}$, and $\phi(y)=\frac{b_1(1+e^{b_1 y})}{1-e^{b_1 y}}$, from which
we obtain the example.

\begin{example}
The Riemannian submersion from Nil space
\begin{align}\notag
\pi  : ( \mathbb{R}^3 , g_{Nil}={\rm d}x^{2}+{\rm
d}y^{2}+({\rm d}z-x{\rm d}y)^{2}) &\to (\mathbb{R}^2 ,dx^2 + (1+x^2)^{-2}dz^2)
\\\notag
\pi(x,y,z) =(x,z)
\end{align}
is not a biharmonic map.
\end{example}
We can check that $e_1=\frac{\partial}{\partial
x},\;e_2=-\frac{x}{\sqrt{1+x^2}}\frac{\partial}{\partial
y}-\sqrt{1+x^2}\frac{\partial}{\partial z}, \;
e_3=\frac{1}{\sqrt{1+x^2}}\frac{\partial}{\partial y}$ form an
orthonormal frame on Nil space adapted to the Riemannian submersion
with $d\pi(e_3)=0, \;\;d\pi(e_i)=\varepsilon_i$, for $i=1,2$ for an
orthonormal frame $\varepsilon_{1}=\frac{\partial}{\partial
x},\;\varepsilon_{2}=-\sqrt{1+x^2}\frac{\partial}{\partial z}$ on
the base space.  We can compute the Lie brackets as
\begin{align*}
[e_1,e_2]=\frac{x}{1+x^2}e_2-\frac{1-x^2}{1+x^2}e_3\\
[e_1,e_3]=-\frac{x}{1+x^2}e_3,\; \;[e_2,e_3]=0,
\end{align*}
from which we obtain the integrability data of the Riemannian
submersion $\pi$ as
$f_1=0,\;f_2=\frac{x}{1+x^2},\;\;\kappa_1=-\frac{x}{1+x^2},\;\;\sigma=\frac{1-x^2}{2(1+x^2)},\;\;\kappa_2=0$.
Since $\kappa_2=0$, we apply Corollary \ref{Cor1} to conclude that
$\pi$ is biharmonic if and only if Equation (\ref{Simp}) holds.
However, a simple computation shows that the left-hand side of the
first equation of (\ref{Simp}) equals $ \frac{x^3-7x}{(1+x^2)^3}$
which does not vanish identically. Thus, the Riemannian submersion
$\pi$ is not biharmonic.

\section{Proofs of the main results}

In this section we will give a complete classification of biharmonic
Riemannian submersions from a $3$-dimensional space form. This is
accomplished by choosing a special adapted orthonormal frame that
exists on a space form which simplifies the biharmonic equation
drastically. We will use the notations $M^3(c)$ for a space form
with constant sectional curvature $c$ and $R_{ijkl}=-\langle
R({e_i},e_j)e_k,e_l\rangle$ for the components of curvature with
respect to an orthonormal basis. The following lemmas will be used
to prove the main theorems.

\begin{lemma}\label{LA}
Let $\pi:M^3(c)\longrightarrow (N^2, h)$ be a Riemannian submersion
from a space form of constant sectional curvature $c$. Then, there
exists an orthonormal frame $\{e_1, e_2, e_3\}$ on $M^3(c)$ adapted
to the Riemannian submersion such that all the integrability data
$f_1, f_2, \kappa_1, \kappa_2$ and $\sigma$ are constant along
fibers of $\pi$, i.e.,
\begin{equation}\label{Hori}
e_3(f_1)=e_3(f_2)=e_3(\kappa_2)= e_3(\kappa_1)= e_3(\sigma)=0.
\end{equation}
\end{lemma}
\begin{proof}
By definition, $f_i=F_i\circ \pi$ for $i=1, 2$, so they are constant
along the fibers. It remains to show that
\begin{equation}\label{RSB66}
e_3(\kappa_2)=0,\;\;\; e_3(\kappa_1)=0,\;\;\; e_3(\sigma)=0.
\end{equation}
One can easily check that the Jacobi identity applies to the frame
$\{e_1, e_2, e_3\}$ yields
\begin{equation}\label{J1}
2e_{3}(\sigma)+\kappa_1f_1+\kappa_2
f_2+e_{2}(\kappa_1)-e_{1}(\kappa_2)=0.
\end{equation}

A straightforward computation using (\ref{J1}) and the fact that
$M^3(c)$ has constant sectional curvature $c$ gives
\begin{equation}\label{RSB67}
\begin{cases}
R^{M}_{1312}=-(e_1(\sigma)-2\kappa_1\sigma)=0,\\
R^{M}_{1313}=-[-e_1(\kappa_1)-\sigma^2+\kappa_{1}^2-\kappa_2f_1]=c,\;\\
R^{M}_{1323}=-[-e_1(\kappa_2)+e_3(\sigma)+\kappa_{1}f_{1}+\kappa_1\kappa_2]=0,\;\\
R^{M}_{1212}=-(e_2(f_1)-e_1(f_2)+f_{1}^{2}+f_{2}^{2}+3\sigma^2)=c,\\
R^{M}_{1223}=-(e_2(\sigma)-2\kappa_1\kappa_2)=0,\\
R^{M}_{2313}=-(-e_2(\kappa_{1})-e_3(\sigma)-\kappa_2 f_2+\kappa_1 \kappa_2)=0,\\
R^{M}_{2323}=-(-\sigma^{2}-e_2(\kappa_2)+\kappa_1f_2+ \kappa_2^2)=c.\\
\end{cases}
\end{equation}

Applying $e_3$ to both sides of the fourth equation of (\ref{RSB67})
and using (\ref{J1}), together with $e_3e_1=[e_3,e_1]+e_1e_3$  and
$e_3e_2=[e_3,e_2]+e_2e_3$, we get
\begin{equation}\notag
\sigma e_3(\sigma)=0,
\end{equation}
which implies
\begin{equation}\notag
e_3(\sigma)=0.
\end{equation}

Using this and applying $e_3$ to both sides of  the 1st and the 5th
equation of (\ref{RSB67}) separately, we obtain
\begin{equation}\notag
e_3(\kappa_{1})=0, \;e_3(\kappa_2)=0,
\end{equation}
which completes the proof of the lemma.
\end{proof}

\begin{lemma}\label{L2}
Let $\pi:(M^3(c),g)\longrightarrow (N^2,h)$ be a Riemannian
submersion with an adapted frame $\{e_1,\; e_2, \;e_3\}$ and the
integrability data $ f_1, f_2, \kappa_1,\;\kappa_2\;{\rm and}\;
\sigma$. Then, there exists another adapted orthonormal frame
$\{e'_1,\;e'_2,e'_3=e_3\}$ on $ M^3(c)$ with integrability data
$f'_1,\;f'_2,\;\kappa'_1=\sqrt{\kappa_{1}^2+\kappa_2^2},\;\kappa_2'
=0, \; {\rm and}\; \sigma'=\sigma$.
\end{lemma}
\begin{proof}
Choose an orthonormal frame $\{e_1,\; e_2, \;e_3\}$ on $M^3(c)$
adapted to the Riemannian submersion $\pi$. It follows from Lemma
\ref{LA} that the integrability data $ \kappa_1 \;{\rm and}\;
\kappa_2$ are constant along the fibers of $\pi$. By a well-known
fact from topology that there exist functions $ \bar{\kappa_1}
\;{\rm and}\; \bar{\kappa_2} \in C^{\infty}(N)$ such that $
\kappa_1= \bar{\kappa_1}\circ \pi\;{\rm and}\;
\kappa_2=\bar{\kappa_2}\circ \pi$. Suppose $e_1, e_2$ are
$\pi$-related to $\varepsilon_1, \varepsilon_2$ respectively. Then,
it is easy to see that
$\varepsilon'_1=\frac{\bar{\kappa_1}}{\sqrt{\bar{\kappa_1}^2+\bar{\kappa_2}^2}}\varepsilon_1
+\frac{\bar{\kappa_2}}{\sqrt{\bar{\kappa_1}^2+\bar{\kappa_2}^2}}\varepsilon_2,
\;\;\varepsilon'_1=\frac{-\bar{\kappa_2}}{\sqrt{\bar{\kappa_1}^2+\bar{\kappa_2}^2}}\varepsilon_1
+\frac{\bar{\kappa_1}}{\sqrt{\bar{\kappa_1}^2+\bar{\kappa_2}^2}}\varepsilon_2$
is an orthonormal frame on the base space. Let $e_1', e_2'$ be the
horizontal lift of $\varepsilon'_1, \varepsilon'_2$ respectively.
Then, one can easily check that the adapted orthonormal frame
$\{e_1',e_2',e_3\}$ satisfies the required conditions stated in the
lemma.
\end{proof}

Now we are ready to give the following classification of biharmonic
Riemannian submersions.
\begin{theorem}\label{MT}
Let $\pi:(M^3(c),g)\longrightarrow (N^2,h)$ be Riemannian submersion
from a space form of constant sectional curvature $c$. Then,
$\pi$ is biharmonic if and only if it is harmonic.\\
\end{theorem}
\begin{proof}
By Lemma \ref{L2}, we can choose an orthonormal frame
$\{e_1,\;e_2,\; e_3\}$ adapted to the Riemannian submersion with
integrability data $\{f_1,\;f_2,\;\kappa_1,\;\kappa_2,\;\sigma\}$
with $\kappa_2=0$. With respect to this frame the curvature equation
(\ref{RSB67}) reduces  to
\begin{equation}\label{RSB67S}
\begin{cases}
e_1(\sigma)-2\kappa_1\sigma=0,\\
-[-e_1(\kappa_1)-\sigma^2+\kappa_{1}^2]=c,\;\\
\kappa_{1}f_{1}=0,\;\\
-(e_2(f_1)-e_1(f_2)+f_{1}^{2}+f_{2}^{2}+3\sigma^2)=c,\\
e_2(\sigma)=0,\\
e_2(\kappa_{1})=0,\\
-(-\sigma^{2}+\kappa_1f_2)=c.\\
\end{cases}
\end{equation}

By the 3rd equation in (\ref{RSB67S}), we have either $\kappa_1=0$
or $f_1=0$. For the first case, $\kappa_1=0$, then, by (\ref{RSB6}),
the tension fields of $\pi$ vanishes and the hence the Riemannian
submersion is harmonic. For the second case, $\kappa_1\ne 0$, then
we have $f_1=0$. We will show that this latter case cannot happen.
We will prove this by using proof by
contradiction in the following two cases:\\

Case I: $\kappa_1\ne 0,\;\; f_1=0$ and $f_2=0$. In this case, the
4th and the 7th equations in (\ref{RSB67S}) implies that
$\sigma=c=0$. Now substituting $f_1=f_2=\sigma=0$ and $\kappa_2=0$
into biharmonic equation (\ref{Simp})  we obtain
\begin{equation}\notag
\Delta\kappa_{1}=0,
\end{equation}
which, by a straightforward computation using (\ref{Conn}), the 2nd,
and the 6th equations of (\ref{RSB67S}), can be turned into
\begin{equation}\notag
\kappa_{1}^3=0.
\end{equation}
It follows that $\kappa_{1}=0$ which is a contradiction.\\

Case II:  $\kappa_1\neq 0,\; f_1=0$ and $f_2\neq 0$. In this case,
we use $f_1=0$ and the 5th, the 6th and the 7th equations of
(\ref{RSB67S}) to reduce the biharmonic equation (\ref{Simp}) into
\begin{equation}\label{d61}
-\Delta^{M}\kappa_1+\kappa_1\{-c-3\sigma^2+f_{2}^{2}\}=0,
\end{equation}
where we have used the fact that the Gauss curvature of the target
surface $K^N=c+3\sigma^2$ obtained from O$'$Neill's curvature
formula for a Riemannian submersion (\cite{On}). A straightforward
computation gives
\begin{eqnarray}\notag
-\Delta^{M}\kappa_1&=&-e_1e_1(\kappa_1)+\nabla_{e_2}e_2(\kappa_1)+\nabla_{e_3}e_3(\kappa_1)\\\notag
&=&-e_1(\kappa_{1}^2-\sigma^2+c)+f_2 e_1(\kappa_1)+\kappa_1
e_1(\kappa_1)\\\notag
&=&5\kappa_1\sigma^2-\kappa_{1}^3-\kappa_{1}c+f_2(\kappa_{1}^2-\sigma^2+c).
\end{eqnarray}
Substituting this into (\ref{d61}) and simplifying the resulting
equation we get
\begin{equation}\label{d81}
\begin{array}{lll}
\kappa_1(3\sigma^2-\kappa_{1}^2-3c)=0.
\end{array}
\end{equation}
By assumption, $\kappa_{1}\neq0$, so (\ref{d81}) implies that
\begin{equation}\label{d9d}
\kappa_{1}^2=3\sigma^2-3c.
\end{equation}
 Applying  $e_1$ to both sides of
 (\ref{d9d}) yields
\begin{equation}\notag
\begin{array}{lll}
\kappa_{1}e_1(\kappa_{1})=3\sigma e_1(\sigma).
\end{array}
\end{equation}
Combining this and the 1st and the 2nd equations in (\ref{RSB67S})
we obtain
\begin{equation}\notag
\begin{array}{lll}
\kappa_{1}(\kappa_{1}^2-\sigma^2+c)=6\kappa_{1}\sigma^2,
\end{array}
\end{equation}
which, since $\kappa_{1}\neq0$, is equivalent to
\begin{equation}\notag
\begin{array}{lll}
(\kappa_{1}^2-\sigma^2+c)=6\sigma^2,
\end{array}
\end{equation}
or
\begin{equation}\label{d131}
\begin{array}{lll}
\kappa_{1}^2=7\sigma^2-c.
\end{array}
\end{equation}
Similarly, applying $e_1$ to both sides of (\ref{d131}) and using
the 1st and the 2nd equations in (\ref{RSB67S}) we get
\begin{equation}\label{d13122}
\begin{array}{lll}
\kappa_{1}^2=15\sigma^2-c.
\end{array}
\end{equation}
Combining (\ref{d9d}), (\ref{d131}) with (\ref{d13122}) we have $
\kappa_{1}=\sigma=c=0$. In particular, $\kappa_{1}=0$, which
contradicts our assumption. Thus, we complete the proof of the
theorem.
\end{proof}

\begin{corollary}
$(1)$ If $\pi:\mathbb{R}^3\longrightarrow (N^2,h)$ is a biharmonic
horizontally homothetic submersion from Euclidean space, then
$(N^2,h)$ is flat and $\pi$ is a composition of an orthogonal
projection $\mathbb{R}^3\longrightarrow \mathbb{R}^2$
followed by a covering map $\mathbb{R}^2\longrightarrow (N^2,h)$;\\
$(2)$ There exists no biharmonic Riemannian submersion
$\pi:\mathbb{H}^3\longrightarrow (N^2,h)$ no matter what $(N^2,h)$
is.
\end{corollary}
\begin{proof}
By a theorem in \cite{OW}, a horizontally homothetic submersion
$\pi:\mathbb{R}^3\longrightarrow (N^2,h)$ is a Riemannian submersion
up to a homothety. It follows from our Theorem \ref{MT} that $\pi$
has to be harmonic and hence a harmonic morphism (see, e.g.,
\cite{BW1}). Using Baird-Wood's Bernstein theorem \cite{BW2} for
harmonic morphisms we conclude that $\pi$ is a composition of an
orthogonal projection $\mathbb{R}^3\longrightarrow \mathbb{R}^2$
followed by a weakly conformal map $\mathbb{R}^2\longrightarrow
(N^2,h)$. Since $\pi$ is a Riemannian submersion, the dilation of
the composition has to be constant $1$, from which we conclude that
the weakly conformal map has to be a covering map and $(N^2,h)$ has
to be flat. This gives the Statement (1). For Statement (2), we
first note that the biharmonic Riemannian submersion
$\pi:\mathbb{H}^3\longrightarrow (N^2,h)$ is a harmonic morphism
because, by Theorem \ref{MT}, it is harmonic. Using again
Baird-Wood's Bernstein theorem \cite{BW3} we conclude that $\pi$ is
the composition of a the orthogonal projection
$\mathbb{H}^3\longrightarrow \mathbb{H}^2$, or the projection to the
plane at infinity $\mathbb{H}^3\longrightarrow \mathbb{C}$, followed
by a weakly conformal map $\rho:\mathbb{H}^2\longrightarrow
(N^2,h)$, or $\rho: \mathbb{C}\longrightarrow (N^2,h)$,
respectively. Suppose the projection has dilation $\lambda_1$ and
the conformal factor of the weakly conformal map is $\lambda_2$,
then, the composition map has dilation $\lambda_1(\lambda_2\circ
\pi)$. It follows from \cite{Gu} that in both cases, $\lambda_1$ is
not constant along the fibers, however, it is clear that
$(\lambda_2\circ \pi)$ is constant along the fibers. It follows that
the product $\lambda_1(\lambda_2\circ \pi)$ cannot be $1$, i.e, in
either case, the map cannot be a Riemannian submersion. This
completes the proof of the corollary.
\end{proof}


\begin{thebibliography}{99}
\bibitem{BW1} P. Baird and  J. C. Wood, {\em Harmonic morphisms between Riemannian manifolds}, London Math. Soc. Monogr.
(N.S.) No. 29, Oxford Univ. Press (2003).
\bibitem{BW2} P. Baird and J. C. Wood, {\em Bernstein
theorems for harmonic morphisms from $R^3$ and $S^3$}, Math. Ann.
280 (1988), no. 4, 579--603.
\bibitem{BW3} P. Baird and J. C. Wood, {\em Harmonic morphisms and conformal foliations by geodesics
of three-dimensional space forms}, J. Austral. Math. Soc. Ser. A 51
(1991), no. 1, 118--153.
\bibitem{CMO2} R. Caddeo, S. Montaldo and C. Oniciuc, {\em Biharmonic submanifolds in spheres},
 Israel J. Math. 130 (2002), 109--123.
\bibitem{Ch} B. Y. Chen, {\em Some open problems and conjectures on submanifolds of finite
type}, Soochow J. Math. 17 (1991), no. 2, 169--188.
\bibitem{CI} B. Y. Chen and S. Ishikawa, {\em Biharmonic pseudo-Riemannian submanifolds
in pseudo-Euclidean spaces}, Kyushu J. Math. 52 (1998), no. 1,
167--185.
\bibitem{Di} I. Dimitri\'c, {\em Submanifolds of $E\sp m$ with harmonic mean curvature
vector}, Bull. Inst. Math. Acad. Sinica 20 (1992), no. 1, 53--65.
\bibitem{Gu} S. Gudmundsson, {\em The geometry of harmonic morphisms},
Ph. D thesis, University of Leeds, 1992.
\bibitem{HV} T. Hasanis and T. Vlachos, {\em Hypersurfaces in $E\sp 4$ with harmonic mean
curvature vector field}, Math. Nachr. 172 (1995), 145--169.
\bibitem{Ji} G. Y. Jiang, {\em $2$-Harmonic maps and their first
and second variational formulas}, Chin. Ann. Math. Ser. A 7(1986)
389-402.
\bibitem{Ji2} G. Y. Jiang, {\em Some non-existence theorems of $2$-harmonic isometric immersions into Euclidean spaces
 }, Chin. Ann. Math. Ser. 8A (1987)
376-383.
\bibitem{LO} E. Loubeau and Y. -L. Ou, {\em Biharmonic maps and
morphisms from conformal mappings}, T$\hat{\rm o}$hoku Math J., to
appear, 2010.
\bibitem{On} B. O'Neill, {\em The fundamental equations of a submersion}, Michigan Math. J. 13 1966
459--469.
\bibitem{OW} Y. -L. Ou and G. Walschap, \emph{A classification of
horizontally homothetic submersions from space forms of nonnegative curvature%
}, Bull. of London Math. Soc., 38(3) (2006), 485-493.
\end{thebibliography}
\end{document}